\documentclass[12pt]{article}
\usepackage{amsmath,amssymb,amsfonts,theorem,makeidx,latexsym,epsfig,,subfigure}

\newtheorem{defn}{Definition}[section]

\newtheorem{lemma}[defn]{Lemma}

{\theorembodyfont{\rmfamily}

\newtheorem{ex}[defn]{Example}}

\newtheorem{thm}[defn]{Theorem}

\newtheorem{prop}[defn]{Proposition}

\newtheorem{cor}[defn]{Corollary}

\numberwithin{equation}{section}

\newcommand{\h}{{\cal H}}

\newcommand{\ltn}{{\ell}^2(\mathbb N)}

\newcommand{\cn}{\mc^n}

\newcommand{\mn}{\mathbb N}

\def\bp{{\noindent\bf Proof. \ }}

\def\ep{\hfill$\square$\par\bigskip}

\def\bqs{\begin{equation}}

\def\eqs{\tag*{$\square$}\end{equation}\par\bigskip}

\def\la{\langle}

\def\ra{\rangle}

\def\ftk{\{f_k\}_{k=1}^\infty}

\def\ctk{\{c_k\}_{k=1}^\infty}

\def\gtk{\{g_k \}_{k=1}^\infty}

\def\htk{\{h_k\}_{k=1}^\infty}

\def\etk{\{e_k\}_{k=1}^\infty}

\def\suk{\sum_{k=1}^\infty}

\def\nl{\left|\left|}

\def\nr{\right|\right|}

\def\span{\overline{\text{span}}}

\def\Span{\text{span}}

\def\vn{\vspace{.1in}\noindent}

\def\bop{\begin{op}\rm}

\def\eop{\end{op}}

\def\cra{{\cal R}}

\def\cn{{\cal N}}

\def\bee{\begin{eqnarray}}

\def\ene{\end{eqnarray}}

\def\bes{\begin{eqnarray*}}

\def\ens{\end{eqnarray*}}

\def\bei{\begin{itemize}}

\def\eni{\end{itemize}}

\def\bt{\begin{thm}}

\def\et{\end{thm}}

\def\bc{\begin{cor}}

\def\ec{\end{cor}}

\def\bpr{\begin{prop}}

\def\epr{\end{prop}}

\def\bl{\begin{lemma}}

\def\el{\end{lemma}}

\def\bd{\begin{defn}}

\def\ed{\end{defn}}

\def\bex{\begin{ex}}

\def\enx{\end{ex}}

\def\bfi{\begin{fig}}

\def\efi{\end{fig}}

\def\eptk{\{\varepsilon_k\}_{k\in\mn}}

\title{Dynamical sampling and frame representations with bounded
operators}

\date{\today}

\author{Ole Christensen, Marzieh Hasannasab, Ehsan Rashidi}

\begin{document}

\maketitle

\begin{abstract} The purpose of this paper is to study frames for a Hilbert space $\h,$
having the form $\{T^n \varphi\}_{n=0}^\infty$ for some $\varphi \in \h$ and an operator
$T:\h \to \h.$ We characterize the frames that have such a representation for a bounded
operator $T,$ and discuss the properties of this operator. In particular, we prove that
the image chain of $T$ has finite length $N$
in the overcomplete case; furthermore $\{T^n \varphi\}_{n=0}^\infty$ has the very particular
property that $\{T^n \varphi\}_{n=0}^{N-1} \cup \{T^n \varphi\}_{n=N+\ell}^\infty$
is a frame for $\h$ for all $\ell\in \mn_0$. We also prove that frames of
the form $\{T^n \varphi\}_{n=0}^\infty$ are sensitive to
the ordering of the elements and to norm-perturbations of
the generator $\varphi$ and the operator $T.$ On the other hand positive stability results are
obtained by considering perturbations of the generator $\varphi$ belonging to an invariant
subspace on which $T$ is a contraction.

\end{abstract}

\section{Introduction}
Let $\h$ denote a separable Hilbert space. A frame is a collection of vectors in $\h$ that
allows each $f\in \h$ to be expanded as an (infinite) linear combination of the
frame elements. Dynamical sampling, as introduced in \cite{A1} by Aldroubi et al., deals with
frame properties of sequences of the form $\{T^n \varphi\}_{n=0}^\infty$,
where $T: \h \to \h$ belongs to certain classes of linear operators and $\varphi \in \h;$ for example, the diagonalizable
normal operators $T$ that lead to a frame for a certain choice of $\varphi$ are characterized.
Further references to dynamical sampling include \cite{A1,A2,A3,CMPP,FP,AK, olemmaarzieh,olemmaarzieh2}.

In this paper we consider the general question of characterizing the frames $\ftk$ for which
a representation of the form $\{T^n \varphi\}_{n=0}^\infty$ with a bounded operator $T: \h \to \h$
exists. While a representation $\{T^n \varphi\}_{n=0}^\infty$ is available for all
linearly independent frames, it is more restrictive to obtain boundedness of
the representing operator.  In Section \ref{70505d} we
give various characterizations of the case where $T$ can be chosen to be bounded: one of them
is in
terms of a certain invariance property of the kernel of the synthesis operator and
another one in terms of a number of equations that must be satisfied.
The result also identifies the unique candidate for the operator $T.$ Several
applications are presented; we prove, e.g., that frames of the form
$\{T^n \varphi\}_{n=0}^\infty$ are sensitive to perturbation of the operator
$T$ and the generator $\varphi,$ and that reorderings of the elements in a frame has a significant influence on the question of being representable by a bounded operator.

The frames having a representation of the form $\{T^n \varphi\}_{n=0}^\infty$ for a bounded
operator naturally split into two classes: the Riesz bases, and certain overcomplete frames.
The operators $T$ appearing for these two classes have very different properties. In the overcomplete case we prove that  there exists some
$N\in \mn_0$ such that $\{T^n \varphi\}_{n=0}^{N-1} \cup \{T^n \varphi\}_{n=N+\ell}^\infty$
is a frame for $\h$ for all $\ell\in \mn_0$.

Section \ref{70928a} collects a number of auxiliary results related to dynamical sampling.
The characterization in Section  \ref{70505d} of operators $T$ for which $\{T^n \varphi\}_{n=0}^\infty$
is a frame identifies $T$ on the form of a mixed frame operator.
We show that there exist frames   $\{T^n \varphi\}_{n=0}^\infty$ for which the operator
$T$ is in fact a frame operator. We also prove that if $\{T^n \varphi\}_{n=0}^\infty$
is  a frame and we perturb
the generator $\varphi$ with an element $\widetilde{\varphi}$ in a $T$-invariant subspace of $\h$ on which $T$ acts as a contraction, then we obtain a
frame $\{T^n (\varphi+ \widetilde{\varphi})\}_{n=0}^\infty$
if the norm of $\widetilde{\varphi}$ is sufficiently small, and otherwise a frame sequence.
Finally, we prove that iterates of
a compact operator acting on a finite collection of vectors can not generate a frame.
This generalizes a result from \cite{A1}.

In the rest of this introduction we will collect the necessary background from frame theory
and operator theory.

\subsection{Frames and operators}
 A sequence $\ftk$
in a Hilbert space $\h$ is a {\it frame}
for  $\h$ if there exist constants $A,B>0$ such that
$$ A \, ||f||^2 \le \suk | \la f, f_k\ra|^2 \le B \, ||f||^2, \, \forall f\in \h.$$
The sequence $\ftk$ is a Bessel sequence if at least the upper frame condition holds.  Also,
it is well-known that for any frame $\ftk$ there exists at least one {\it dual frame}, i.e.,
a frame $\gtk$ such that
\bes f= \suk \la f,g_k \ra f_k, \, \forall f\in \h.\ens

A sequence $\ftk$ in $\h$ is a {\it Riesz basis} if $\span \ftk=\h$ and
there exist constants $A,B>0$ such that
$$A \suk |c_k|^2 \le \nl \suk c_k f_k \nr^2 \le B\suk |c_k|^2$$
for all finite scalar sequences $\ctk,$ i.e., sequences where only
a finite number of entries are nonzero. A Riesz basis is automatically a frame; and
a frame is a Riesz basis if and only if it is $\omega$-independent, meaning that $\suk c_kf_k$ only
vanishes if $c_k=0, \, \forall k\in \mn.$

Throughout the paper we also need to consider {\it linearly independent frames;} here linear
independence should be understood in the classical linear algebra sense, i.e., that a finite
linear combination of frame elements only vanishes if all coefficients vanish.

If $\ftk$ is a Bessel sequence, the {\it synthesis operator} is defined by
\bee \label{60811a} U: \ltn \to \h, \, U\ctk := \suk c_k f_k;\ene
it is well known that $U$ is well-defined and bounded.
A central role will be played by the kernel
of the operator $U,$ i.e., the subset of $\ltn$ given by
\bee \label{60811f} \cn(U)=\left\{\ctk \in \ltn ~\bigg|~\suk c_kf_k=0\right\}.\ene

Let us state  an example of a frame that indeed has the form
$\{ T^n \varphi\}_{n=0}^\infty$ for a bounded operator $T: \h\to \h.$ We will refer
to this example at several places in the sequel.
\bex \label{70614b} Based on interpolation sequences in the Hardy space on the
unit disc, Aldroubi et al. \cite{A2} have constructed frames $\{ T^n \varphi\}_{n=0}^\infty$
for $\ell^2(\mn).$ We will formulate the result in the setting of a separable
infinite-dimensional Hilbert space $\h$. Consider an operator $T$ of the form $T= \sum_{k=1}^\infty \lambda_k P_k,$ where
$P_k, k\in \mn,$ are rank 1 orthogonal projections such that $P_jP_k=0, \, j\neq k, \, \suk P_k=I,$
and $|\lambda_k| <1$ for all $k\in \mn.$ Then
there exist unit vectors $e_k$ such that $P_kf=\la f, e_k\ra e_k.$ The condition $P_jP_k=0$
implies that $\la e_j, e_k\ra=0$ for $j\neq k;$
since $\suk P_k=I,$ it follows that $\{e_k\}_{k=1}^\infty$ is an orthonormal
basis for $\h.$
Assume that $\{\lambda_k\}_{k=1}^\infty$ satisfies the so-called Carleson condition, i.e.,
\bee\label{carleson}\displaystyle\inf_k \prod_{j\neq k}\frac{|\lambda_j-\lambda_k|}{|1-\lambda_j\overline{\lambda_k}|}>0.\ene
Then, letting
$\varphi:= \sum_{k=1}^\infty \sqrt{1-|\lambda_k|^2} e_k,$ the family $\{T^n\varphi\}_{n=0}^\infty $ is a frame for $\h.$
A concrete example of a sequence satisfying the Carleson condition is $\{\lambda_k\}_{k=1}^\infty=\{1-\alpha^{-k}\}_{k=1}^\infty$ for some $\alpha >1.$
\ep \enx
For more information about frames we refer to the monographs \cite{Heil,CB}.

Finally, we will need the right-shift operator on $\ltn,$ defined by
\bee \label{70606a} {\cal T}: \ltn \to \ltn, {\cal T}\ctk= \{0, c_1, c_2, \dots\}.\ene

For any operator $T: \h \to \h$, the null spaces of $T^n, \, n\in \mn_0,$ form an increasing sequence $\cn(T^0)\subseteq\cn(T)\subseteq \cn(T^2)\subseteq\cdots$. We say that the {\it null chain} is finite if there is an $N\in\mn_0$ such that $\cn(T^N)=\cn(T^{N+1})$; in that case the smallest
 such $N$ is called the {\it length} of the null chain. The {\it image chain} of $T$ is the decreasing sequence $\cra(T^0)\supseteq\cra(T^1)\supseteq\cra(T^2)\supseteq\cdots$. The image chain of $T$ is finite if there is an $N\in\mn_0$ such that $\cra(T^N)=\cra(T^{N+1})$; in that case the smallest
 such $N$ is called the {\it length} of the image chain, and is denoted by $q(T).$

 \section{Boundedness of the operator $T$} \label{70505d}
In the entire section $\h$ denotes an infinite-dimensional separable Hilbert space. In
\cite{olemmaarzieh} it is proved that a frame  $\ftk$ for $\h$ has a representation
$\ftk= \{T^nf_1\}_{n=0}^\infty$ for some  operator $T: \h \to \h$ if and only
if $\ftk$ is linearly independent. In the following Theorem \ref{70510a} we give  various characterizations
of the case where such a representation is possible with a bounded operator $T$.
Note that in the affirmative case, the operator $T$ is unique; the result also  identifies the only possible candidate for the operator $T$, see
\eqref{70505b}.
The flavour of the characterizations in Theorem \ref{70510a} are quite different. Indeed,
the characterization in  (iii) deals with an invariance
property  of the kernel of the synthesis operator. On the other hand, (ii) is an ``intrinsic characterization" that is formulated directly in terms of the elements in the given frame $\ftk:$
if we have access to a dual frame, it allows to check the existence of a bounded operator
$T$ by verifying a number of equations. We illustrate both characterizations in Example \ref{70513p}.

		\bt\label{70510a}
Consider a  frame  $\ftk$ with frame bounds $A,B.$ Then the following are equivalent:

\bei
\item[(i)] The frame has a representation $\ftk=\{T^nf_1\}_{n=0}^\infty$ for some
bounded operator $T:\h \to \h.$
\item[(ii)] For some dual frame $\gtk$ (and hence all),
\bee \label{70510b} f_{j+1}= \suk \la f_j, g_k\ra f_{k+1}, \, \forall j\in \mn.\ene

\item[(iii)] The kernel $\cn(U)$
  of the synthesis operator $U$ is invariant  under the right-shift operator ${\cal T}.$
  \eni In the affirmative case, letting $\gtk$ denote an arbitrary dual frame of $\ftk,$
the operator $T$ has the form
\bee \label{70505b}  Tf= \sum_{k=1}^\infty \la f, g_k\ra f_{k+1}, \, \forall f\in \h,\ene
and $1 \le \| T\|\le \sqrt{BA^{-1}}.$\et

\bp We first note that the only possible candidate for a bounded operator $T$
representing the frame $\ftk$ indeed
is the one given in \eqref{70505b}; in fact, if $\ftk=\{T^nf_1\}_{n=0}^\infty$ and $T$
is bounded, just apply the operator $T$ on the decomposition $f=\sum_{k=1}^\infty \la f, g_k\ra f_k.$

\noindent (i)$\Rightarrow$ (ii). This follows directly from the observation above
by applying \eqref{70505b} on $f=f_j.$

\noindent (ii)$\Rightarrow$ (i).
The only possible choice of a representing operator
$T$ is given in \eqref{70505b}, and \eqref{70510b} is expressing precisely
that $Tf_j=f_{j+1}$ for all $j\in \mn,$ i.e., that $\ftk = \{T^nf_1\}_{n=0}^\infty.$

\noindent (i)$\Rightarrow$ (iii). Assume that
 $\ctk\in \cn(U)$. Using the boundedness of $T$, a direct calculation shows that $U\mathcal{T}\ctk=TU\ctk=T0=0$, which means that $\mathcal{T}\ctk\in \cn(U).$
 				
\noindent (iii) $\Rightarrow$ (i). We will first show that the condition (iii) implies
that $\ftk$ is linearly independent. Assume that $\sum_{k=1}^Nc_kf_k=0$ for some
$N\in \mn,$ and that
$c_N\neq 0.$ Then $f_N \in \Span\{f_1, \dots, f_{N-1}\}.$ Then, using the ${\cal T}$-invariance
of $\cn(U),$ we have that $\sum_{k=1}^N c_kf_{k+1}=0,$ which implies
that $f_{N+1}\in \Span\{f_1, \dots, f_{N}\}\\= \Span\{f_1, \dots, f_{N-1}\}.$ By induction,
it follows that $f_{N+\ell}\in \Span\{f_1, \dots, f_{N}\}$ for all $\ell \in \mn,$ but
this contradicts that $\ftk$ is a frame for an infinite-dimensional space. Thus $c_N=0,$
which by iteration implies that $0=c_N=c_{N-1}= \cdots= c_1.$ Thus $\ftk$ is linearly
independent, as claimed.

The linear independence of $\ftk$ implies that we can define a linear operator $T:\Span\ftk\to\Span\ftk$ by $Tf_k=f_{k+1}$. We now
 show that $T$ is bounded if  $\cn(U)$ is invariant  under right-shifts. Let $f=\sum_{k=1}^{N}c_kf_k$ for some $\ctk\in \ltn$ with $c_k=0$ for $k\geq N+1$. Let us decompose $\ctk$ as  $\ctk=\{d_k\}_{k=1}^{\infty}+\{r_k\}_{k=1}^{\infty}$, with $\{d_k\}_{k=1}^{\infty}\in \cn(U)$ and $\{r_k\}_{k=1}^{\infty}\in \cn(U)^{\perp}$. Since $\cn(U)$ is invariant under right-shifts, we have
				$\sum_{k=1}^{\infty}d_kf_{k+1}=0$. Since $\ftk$ is a Bessel sequence, it follows that
				\bes \|Tf\|^2=\big\|\sum_{k=1}^{N}c_kf_{k+1}\big\|^2=
\big\|\sum_{k=1}^{\infty}r_kf_{k+1}\big\|^2
				\leq  B\sum_{k=1}^{\infty}|r_k|^2. \ens
Using Lemma 5.5.5 in \cite{CB}, since $\{r_k\}_{k=1}^{\infty}\in \cn(U)^{\perp}$ and
 $\ftk$  is a frame with lower bound
$A,$ we have
$A \sum_{k=1}^{\infty} |r_k|^2 \le \|\sum_{k=1}^{\infty}r_kf_{k}\|^2.$ It follows that
				\bes
				\left\|Tf\right\|^2
			&	\leq& BA^{-1} \big\|\sum_{k=1}^{\infty}r_kf_{k}\big\|^2 =BA^{-1}\big\|\sum_{k=1}^{\infty}(d_k+r_k)f_{k}\big\|^2 \\
			& =&  BA^{-1} \big\|\sum_{k=1}^{\infty}c_kf_{k}\big\|^2
				=  BA^{-1}\left\|f\right\|^2.
				\ens
				Therefore $T$ can be extended to a bounded operator on $\h$, as claimed.
The proof shows that $||T|| \le \sqrt{BA^{-1}};$ the inequality $||T||\ge 1$ is proved in
\cite{A3}.
\ep

Theorem \ref{70510a} gives an easy proof of a result from \cite{olemmaarzieh} :

\bex \label{70510d} Let $\ftk$ be a Riesz basis for $\h,$ with dual Riesz basis $\gtk$.
Then the condition \eqref{70510b} is trivially satisfied. Thus $\ftk$ has
a representation on the form  $\{T^n f_1\}_{n=0}^\infty$
for a bounded operator $T.$
\ep \enx

In the next example we give one more application of Theorem \ref{70510a}.
We show that (i): the
set of bounded operators $T: \h \to \h$ for which $\{T^n \varphi\}_{n=0}^\infty $ is a frame for
some $\varphi \in \h$ does not form an open set in $B(\h),$ and (ii): the set of $\varphi \in \h$ for which
$\{T^n \varphi\}_{n=0}^\infty $ is a frame for a fixed bounded operator $T$ does not form an open set in $\h.$ In a certain sense these results explain the fact that the known frames of the form
$\{T^n \varphi\}_{n=0}^\infty $ with a bounded operator $T$ are very particular.

\bex \label{70929a}  (i) Consider any frame $\{T^n \varphi\}_{n=0}^\infty$ for which $T:\h \to \h$ is bounded
and $||T||=1.$ Now, given any $\epsilon\in (0,1),$ define the operator $W: \h \to \h$ by
$W:= (1-\epsilon)T.$ Then $||T-W||= || \epsilon T||= \epsilon.$  However,
$||W|| = 1-\epsilon<1,$ which by Theorem \ref{70510a} implies that
$\{W^n \psi\}_{n=0}^\infty$ is not a frame for any $\psi \in \h.$ This proves that
the set
\bes \{T\in B(\h) \, \big| \, \mbox{there exists }\varphi\in \h \mbox{ such that } \{T^n \varphi\}_{n=0}^\infty \mbox{ is a frame } \} \ens
is not open in $B(\h).$

\vn (ii) Given a bounded operator $T: \h \to \h$, the set
$\{ \varphi \in \h \, \big| \, \{T^n \varphi\}_{n=0}^\infty \mbox{   is a frame  }\}$
is not open in general. In order to illustrate this, consider the frame $\{T^n\varphi\}_{n=0}^\infty $ in Example \ref{70614b}.
Given any $\epsilon >0,$ choose $\ell\in \mn$ such that $\sqrt{1-|\lambda_\ell|^2} < \epsilon,$
and let $\psi:= \sum_{k\neq \ell} \sqrt{1-|\lambda_k|^2} e_k= \varphi - \sqrt{1-|\lambda_\ell|^2} e_\ell.$
Then $|| \varphi - \psi|| < \epsilon,$ but
$\{T^n\psi\}_{n=0}^\infty $ is not a frame. Thus, the set of functions $\varphi$ generating
a frame for a fixed operator $T$ is not open. \ep \enx

In general, the operator $T$ in \eqref{70505b}  depends on the choice
of the dual frame $\gtk;$ however, in the case where the frame $\ftk$ indeed has the desired
representation in terms of a bounded operator, the operators in \eqref{70505b}
become independent of the choice of $\gtk.$ In other words:
we can falsify the existence of a bounded operator representing the frame by calculating
the operator in \eqref{70505b} for two different dual frames, and show that they are not
equal.

\bex Let $\etk$ denote an orthonormal basis for $\h$ and consider the frame
$\ftk= \{e_1, e_1,e_2, e_3, \dots\}.$ Considering the dual frame $\{0, e_1, e_2, e_3, \dots\},
$ the operator $T$ in \eqref{70505b} takes the form $T_1f= \suk \la f, e_k\ra e_{k+1};$ for
the dual frame $\{e_1, 0,  e_2, e_3, \dots\},
$ we obtain  $T_2f= \la f, e_1\ra e_1+\sum_{k=2}^\infty \la f, e_k\ra e_{k+1}.$ Since
$T_1\neq T_2,$ this shows that $\ftk$ does not have a representation in terms
of a bounded operator. \ep \enx

We will give another application of Theorem \ref{70510a} in Example \ref{70513p}, but
let us first state a rather surprising result about frames of the
form $\{T^n \varphi\}_{n=0}^\infty$.

It is well-known that if a family $\ftk$ is a frame for $\h,$ then the subfamily $\{f_k\}_{k=N}^\infty$
is a frame for the subspace $V_N:=\span \{f_k\}_{k=N}^\infty$ for all $N\in \mn.$
When we increase $N$ it corresponds to remove more elements from the original frame $\ftk,$ so in general
we expect the spaces $V_N$ to become smaller. However, for overcomplete frames having a representation of the form
$\{T^n \varphi\}_{n=0}^\infty$ with a bounded operator $T,$ the spaces $V_k$ stabilizes at some
point, and the sequence $\{T^n \varphi\}_{n=N}^\infty$ remains a frame for the same space, no
matter how large $N$ is chosen to be. Intuitively, this can be formulated by saying that
the frame $\ftk$ has infinite excess in ``almost all directions".
The result is based on the proof of the following proposition which is of independent interest.

\bpr \label{70428b} Assume that $\{f_k\}_{k=1}^\infty$ is an overcomplete frame and that
$\{f_k\}_{k=1}^\infty= \{T^n\varphi\}_{n=0}^\infty$ for some bounded linear
operator $T: \h \to \h.$ Then there exists an $N\in \mn$ such that
\bee \label{70826a} \{f_k\}_{k=1}^N \cup \{f_k\}_{k=N+\ell}^\infty
\ene is a frame for $\h$ for all $\ell\in \mn.$ \epr

\bp Choose
some coefficients $\{c_k\}_{k=1}^\infty\in \ell^2(\mn)$ such that
$\sum_{k=1}^\infty c_kf_k=0.$ Letting $N:= \min\{k\in \mn\, | \, c_k \neq 0\},$
we have that
\bee \label{70428a} -c_Nf_N= \sum_{k=N+1}^\infty c_kf_k,\ene so
$f_N\in \span \{f_k\}_{k=N+1}^\infty.$ Thus $\{f_k\}_{k=N+1}^\infty$
is a frame for $\span \{f_k\}_{k=N}^\infty.$ Applying
the operator $T$ on \eqref{70428a} shows that $f_{N+1}\in \span \{f_k\}_{k=N+2}^\infty.$ By
iterated application of the operator $T$ this proves that for any $\ell\in \mn,$ the family $\{f_k\}_{k=N+\ell}^\infty$ is
a frame for $\span \{f_k\}_{k=N}^\infty,$ which leads to the desired result. \ep

\begin{thm}\label{2306a} Assume that $\{T^n \varphi\}_{n=0}^\infty$ is an overcomplete frame for some $\varphi\in\h$
and some bounded operator $T: \h \to \h.$ Then the following hold:
\bei
\item[(i)]  The image chain for the operator $T$ has finite length $q(T).$
\item[(ii)] If $N\in \mn_0,$ then
$T^N\varphi\in\span\{T^n\varphi\}_{n=N+1}^\infty \Leftrightarrow N\ge q(T)$.
\eni
For any $N\ge q(T),$ let $V:=\span\{T^n \varphi\}_{n=N}^\infty$. Then the following hold:
\bei
\item[(iii)] The space $V$ is independent
of $N$ and has finite codimension.
\item[(iv)] The sequence $\{T^n \varphi\}_{n=N+ \ell}^\infty$ is a frame for $V$
for all $\ell\in \mn_0.$
\item[(v)] $V$ is invariant under $T$, and $T:V\to V$ is surjective.
\item[(vi)]  If the null chain of $T$ has finite length then $T:V\to V$ is injective; in particular this is the case if  $T$ is normal.
 \eni  \end{thm}

\bp
$(i)$ For $k\in \mn_0,$ let $V_k:= \span \{T^n \varphi\}_{n=k}^\infty.$ Then $V_{k+1} \subseteq V_{k}$
for all $k\in \mn_0.$
By the proof of Proposition \ref{70428b}  there exists an $N\in \mn_0$ such that
$V_N= V_{N+\ell}$ for all $\ell\in \mn.$
Since
\bes \cra(T^k)= \left\{T^kf \, \big| \, f\in \h \right\} & = & \left\{T^k\sum_{n=0}^\infty c_n T^n \varphi \, \big|
\{c_n\}_{n=0}^\infty \in \ell^2(\mn_0)\right\} \\ & = & \left\{\sum_{n=0}^\infty c_n T^{k+n} \varphi \, \big|
\{c_n\}_{n=0}^\infty \in \ell^2(\mn_0)\right\} =V_k,\ens it follows immediately  that the image chain of $T$ has finite length. \\
$(ii)$  If $T^N \varphi\in\span\{T^n \varphi\}_{n=N+1}^\infty=V_{N+1}$ for some $N\in \mn_0,$ then $V_N= V_{N+1}$. Similarly to the proof of Proposition \ref{70428b}  it follows
that $V_N=V_{N+\ell}$ for all $\ell \in \mn_0$, and hence
$\cra(T^N)= \cra(T^{N+\ell})$ for all $\ell \in \mn_0;$ in particular, $N\geq q(T)$.
On the other hand, if $N\geq q(T)$, then the fact that
$\cra(T^k)= V_k$ implies that $V_N=V_{N+1};$ thus  $T^N \varphi\in\span\{T^n \varphi\}_{n=N+1}^\infty$, as desired.\\
$(iii)$ That the space $V$ is independent of $N$ as long as $N\ge q(T)$ was
proved in (ii); furthermore the definition of $V$ shows that it has finite codimension.\\
(iv) We already saw this in the proof of Proposition \ref{70428b}.  \\
$(v)$  The definition of $V$ shows that it  is invariant under $T.$ Using that $V=V_{N+1}$, we can write any $f\in V$ on the form
$$f= \sum_{n=N+1}^\infty c_n T^n \varphi= \sum_{n=N}^\infty c_{n+1} T^{n+1} \varphi
=T\sum_{n=N}^\infty c_{n+1} T^{n} \varphi$$ for some $\ctk\in \ltn;$ thus $T:V \to V$ is
surjective.\\
$(vi)$ By assumption the null chain of $T$ has finite length; since the image chain also has finite length by $(i)$, the lengths are equal by Proposition 3.8 in \cite{Heu}. It follows that $\cn(T^N)=\cn(T^{N+1})$ whenever $N \ge q(T)$. Now assume that $Tf=0$ for some $f\in V$. Since $\{T^{N+n}\varphi\}_{n=0}^\infty$ is a frame for $V$, there exists a sequence  $\{c_n\}_{n=0}^\infty\in\ell^2(\mn_0)$ such that $f=\sum_{n=0}^{\infty}c_n T^{N+n}\varphi$. Therefore $T^{N+1}(\sum_{n=0}^{\infty}c_n T^{n}\varphi)=Tf=0$ which means that $\sum_{n=0}^{\infty}c_n T^{n}\varphi\in \cn(T^{N+1})=\cn(T^N)$. This implies that $$f=\sum_{n=0}^{\infty}c_n T^{N+n}\varphi=T^N\sum_{n=0}^{\infty}c_n T^{n}\varphi=0.$$
Thus $T$ is injective, as claimed.
 Assuming now that $T$ is normal, the sequence $\{(T^*)^n\varphi\}_{n=0}^\infty$ is also a frame
by \cite{CMPP} and since $T^*$ is normal, it is overcomplete by \cite{A2}. Applying the result in $(i)$ to $T^*$ shows that there is some $M\in\mn$ such that $\cra((T^*)^M)=\cra((T^*)^{M+1})$.  Therefore  $\cn(T^M)=\cn(T^{M+1})$. By the first part of $(v)$ it now follows that $T$ is injective on $V$.
\ep

The following example shows a case where the assumptions in Theorem \ref{2306a} do not imply that $T$ is injective considered as an operator on $\h$, even though it is bijective on the invariant subspace $V$.
\bex
Let us return to the setup in Example \ref{70614b}. If a sequence $\{\lambda_k\}_{k=1}^\infty$
in the unit disc
satisfies the Carleson condition and consists of nonzero numbers, then also
$\{0\} \cup \{\lambda_k\}_{k=1}^\infty$ satisfies the Carleson condition. Thus,
without loss of generality we can assume that there is some $K\in\mn$ such that $\lambda_K=0$. Then clearly the operator $T$ is not injective on $\h$ but it is bijective on the subspace $V=\span\{T^n\varphi\}_{n=K+1}^\infty$.
\ep\enx

The following example illustrates Theorem \ref{70510a} and also provides insight
in Theorem \ref{2306a}. We will state a number of consequences after the example itself.

\bex \label{70513p} Let $\h_1$ denote a finite-dimensional Hilbert space and $\h_2$ an infinite-dimensional
separable Hilbert space. Let $\{e_k\}_{k=1}^N$ denote a (Riesz) basis for $\h_1,$ let
$\htk$ be a frame for $\h_2$ , and consider the sequence $\ftk$
in $\h:= \h_1 \oplus \h_2$ given by
\bes \ftk= \{e_1, e_2, \dots, e_N, h_1, h_2, \dots\}.\ens
Assuming that $\htk$ has a representation in terms of a bounded operator
as in Theorem \ref{70510a}, we want to show that there exists a bounded
operator $T: \h \to \h$ such that $\ftk=\{T^nf_1\}_{n=0}^\infty.$

Let us first do so by verifying the condition \eqref{70510b} in Theorem \ref{70510a}.
Letting $\{\widetilde{e}_k\}_{k=1}^N$ denote the dual Riesz basis for $\{e_k\}_{k=1}^N$ and $\{\widetilde{h}_k\}_{k=1}^\infty$  a dual frame for $\htk$, the sequence
$\gtk=\{\widetilde{e}_1,\widetilde{e}_2,\dots,\widetilde{e}_N,
\widetilde{h}_1,\widetilde{h}_2,\cdots \}$
is a dual frame for the frame $\ftk$. Now, for $j=1, \dots,N,$
using that  $\{e_k\}_{k=1}^N$ and $\{\widetilde{e}_k\}_{k=1}^N$ are biorthogonal
and that $\h_1 \bot \h_2,$ it follows that \eqref{70510b} holds for $j=1,\cdots,N$. For $j>N$, we have $f_j=h_{j-N}$, which  is perpendicular to the first $N$ elements of $\gtk.$
Using Theorem \ref{70510a} (ii) on the frame $\htk,$ it follows that
\bes\label{1705a} \sum_{k=1}^{\infty}\langle f_j, g_k\rangle f_{k+1}=\sum_{k=1}^{\infty}\langle h_{j-N}, \widetilde{h}_k\rangle h_{k+1}=h_{j-N+1}=f_{j+1}; \ens
thus we have also verified \eqref{70510b} for $j>N$, as desired.

Let us give an alternative proof using the condition in Theorem \ref{70510a} (iii).
We first note that since $\{e_k\}_{k=1}^N$
and $\htk$ are linearly independent sequences in orthogonal spaces, $\ftk$ is linearly
independent, and hence representable on the form $\ftk= \{T^n f_1\}_{n=0}^\infty$
for a linear operator $T: \Span \ftk \to \h.$ In order to show
that $T$ can be
extended to a bounded operator on $\h,$
consider a sequence $\ctk\in \cn(U)$; then
\[ 0=\suk c_k f_k =\sum_{k=1}^{N} c_k e_k + \sum_{k=N+1}^\infty c_{k}h_{k-N}. \]
Since $\h_1 \bot \h_2,$ this implies that
$\sum_{k=1}^{N} c_k e_k =\sum_{k=N+1}^\infty c_{k}h_{k-N}=0.$
It immediately follows that $c_1=c_2=\cdots=c_N=0$. Also, applying   Theorem \ref{70510a} (iii)
on the sequence $\htk,$ we conclude that $\sum_{k=N+1}^\infty c_k h_{k-N+1}=0$, i.e., $\sum_{k=N+1}^\infty c_k f_{k+1}=0$. It follows that $\sum_{k=1}^{\infty}c_k f_{k+1}=0,$
i.e., that  $\cn(U)$ is invariant under right-shifts, as desired.
\ep \enx

The construction in Example \ref{70513p} is useful for various purposes.
\bei \item
Example \ref{70513p} illustrates that in the setup of Proposition \ref{70428b} it is necessary to include
the vectors $\{f_k\}_{k=1}^N$ in \eqref{70826a}; indeed, letting $\htk$ in
Example  \ref{70513p} be a overcomplete frame for $\h_2,$ the frame
$\ftk$ is overcomplete in $\h,$ but the vectors $\{e_k\}_{k=1}^N$ are not redundant
and can not be removed if we want to keep the frame property.
\item Example \ref{70513p} shows that there exist overcomplete frames $\ftk$
that are represented by an operator $T$ for which $||T||>1;$ such a construction
is obtained, e.g., by letting $\htk$ in
Example  \ref{70513p} be a overcomplete frame for $\h_2,$ and choosing the Riesz
basis $\{e_k\}_{k=1}^N$ such that it is represented by an operator
with norm strictly larger than one.
\eni

We will give yet another application of Theorem \ref{70510a}, showing that the
ordering of the elements in a frame is very important for the question of being representable
by a bounded operator. Indeed, if $\ftk$ has such a representation, a simple reordering
of  two elements $f_\ell$ and $f_{\ell^\prime}$ destroys this property if
$\span \{f_k\}_{k\notin\{\ell-1,
\ell, \ell^\prime-1,
\ell^\prime\}}=\h$. Notice that by Theorem \ref{2306a} this assumption on $\ell, \ell^\prime$ is very weak and only excludes a finite number of choices
of the elements $f_\ell$ and $f_{\ell^\prime}$. However, Example \ref{70513p} also
demonstrates that the assumption is necessary, in the sense that
the result does not hold in general: in that example we can clearly change
the order of any of the vectors $\{e_1, \dots, e_N\}$ without affecting the boundedness of
the representing operator.

\bc Assume that the frame $\ftk$ has a representation $\{T^nf_1\}_{n=0}^\infty,$ where $T: \h \to \h$ is bounded. Choose $\ell \neq \ell^\prime$ such that $\span \{f_k\}_{k\notin\{\ell-1,
\ell, \ell^\prime-1,
\ell^\prime\}}=\h$, and let $\{\widetilde{f_k}\}_{k=1}^\infty$ denote the sequence consisting
of the same elements as $\ftk$ but with $f_\ell$ and $f_{\ell^\prime}$ interchanged. Then
$\{\widetilde{f_k}\}_{k=1}^\infty$ is not representable by a bounded operator. \ec

\bp Let us apply the characterization of boundedness in Theorem \ref{70510a} (ii) with
$\gtk$ being the canonical dual frame of $\ftk.$ Then, by the assumption on the frame $\ftk,$
\bee \label{70630a} f_{j+1} & = & \sum_{k=1}^\infty \la f_j, g_k\ra f_{k+1} = \sum_{k\notin\{\ell-1,
\ell, \ell^\prime-1,
\ell^\prime\}} \la f_j, g_k\ra f_{k+1} \\ \notag & \ & +\la f_j, g_{\ell-1}\ra f_{\ell}+ \la f_j, g_\ell\ra f_{\ell+1}+\la f_j, g_{\ell^\prime-1}\ra f_{\ell^\prime}+
\la f_j, g_{\ell^\prime}\ra f_{\ell^\prime+1}.\ene
Now consider the sequence  $\{\widetilde{f_k}\}_{k=1}^\infty$; its canonical dual frame
$\{\widetilde{g_k}\}_{k=1}^\infty$
consists of the same elements as $\gtk,$ but with the order of $g_\ell$ and $g_{\ell^\prime}$
interchanged. In order
to reach a contradiction, assume that  $\{\widetilde{f_k}\}_{k=1}^\infty$ is also
representable by a bounded operator. Then, for $j\notin \{\ell-1,
\ell, \ell^\prime-1,
\ell^\prime\},$ a new application of Theorem \ref{70510a} (ii) yields that
\bee \label{70630b} f_{j+1}= \widetilde{f_{j+1}} & = &  \sum_{k=1}^\infty \la  \widetilde{f_j},  \widetilde{g_k}\ra \widetilde{f_{k+1}} = \sum_{k\notin\{\ell-1,
\ell, \ell^\prime-1,
\ell^\prime\}} \la f_j, g_k\ra f_{k+1} \\ \notag & \ & +\la f_j, g_{\ell-1}\ra f_{\ell^\prime}+ \la f_j, g_{\ell^\prime}\ra f_{\ell+1}+\la f_j, g_{\ell^\prime-1}\ra f_{\ell}+
\la f_j, g_{\ell}\ra f_{\ell^\prime+1}.\ene Comparing \eqref{70630a} and \eqref{70630b}
shows that for $j\notin \{\ell-1,
\ell, \ell^\prime-1,
\ell^\prime\},$
\bes & \ & \la f_j, g_{\ell-1}\ra f_{\ell}+ \la f_j, g_\ell\ra f_{\ell+1}+\la f_j, g_{\ell^\prime-1}\ra f_{\ell^\prime}+
\la f_j, g_{\ell^\prime}\ra f_{\ell^\prime+1} \\ & = &
\la f_j, g_{\ell-1}\ra f_{\ell^\prime}+ \la f_j, g_{\ell^\prime}\ra f_{\ell+1}+\la f_j, g_{\ell^\prime-1}\ra f_{\ell}+
\la f_j, g_{\ell}\ra f_{\ell^\prime+1}.
 \ens Without loss of generality, assume that $\ell> \ell^\prime;$
 by the linear independence of the elements in $\ftk$ this in particular implies that
$\la f_j, g_\ell\ra = \la f_j, g_{\ell^\prime}\ra.$ Using that
$\span \{f_k\}_{k\notin\{\ell-1,
\ell, \ell^\prime-1,
\ell^\prime\}}=\h$, we conclude that $g_\ell= g_{\ell^\prime}$ and therefore,
by applying the frame operator, $f_\ell= f_{\ell^\prime}.$ However, this contradicts the
linear independence of the elements in $\ftk;$ thus, we conclude that
$\{\widetilde{f_k}\}_{k=1}^\infty$ is not representable by a bounded operator.\ep

\section{Auxiliary results} \label{70928a}
In this short section we provide a few more results  related to the results in Section
\ref{70505d}.

\subsection{Tight frames}

\bc Consider a frame $\{T^n\varphi\}_{n=0}^\infty,$ where $T: \h \to \h$ is bounded.
Then the following hold:
\bei
\item[(i)] If $\{T^n\varphi\}_{n=0}^\infty$ is a tight frame, then $||T||=1.$
\item[(ii)] The canonical tight frame associated with $\{T^n\varphi\}_{n=0}^\infty$ is
$\{ S^{-1/2}T^n\varphi\}_{n=0}^\infty  =
 \{ \left(S^{-1/2}T S^{1/2}\right)^n S^{-1/2}\varphi\}_{n=0}^\infty,$
where $S:\h \to \h$ is the frame operator; in particular, $\nl S^{-1/2}T S^{1/2}\nr =1.$
\item[(iii)] If $||Tf||=c\, ||f||$ for all $f\in \h,$ then $c=1,$ i.e., $T$
is isometric.
\eni \ec
\bp (i) and the fact that $c\ge1$ in (iii) follow immediately from the norm-estimate in Theorem \ref{70510a}; the proof of (iii) is completeted by noticing that if $c>1,$ then
$||T^n\varphi|| = c^n\, ||\varphi||\to \infty$ as $n\to \infty,$ which violates the frame
property. The result in (ii) follows by direct calculation and an application of (i).\ep

\subsection{Iterated systems and the frame operator}

Returning to Theorem \ref{70510a}, we note
that \eqref{70505b} identifies the only possible candidate for an operator $T$
representing a frame $\ftk$ in form of a mixed frame operator. Certain frames $\ftk$
are indeed represented by a frame operator:

\bex\label{6057a} In Example \ref{70614b} we considered frames $\{T^n \varphi\}_{n=0}^\infty$
for operators of the form $T=\sum \lambda_k P_k,$ where $P_k$ are rank 1 projections such
that $P_jP_k=0$ for $j\neq k$ and $\suk P_k=I.$ Choosing the orthonormal basis $\etk$
such that $P_kf= \la f,e_k\ra e_k,$ and considering the case where $\lambda_k\ge 0$ for all
$k\in \mn,$ it is clear that $T$ is the frame operator for the sequence $\{\sqrt{\lambda_k} e_k \}_{k=1}^\infty.$
\ep \enx

In general, a bounded operator $T: \h \to \h$ is a frame operator if and only if
$T$ is positive and invertible.

\subsection{Perturbation of a frame $\{T^n \varphi\}_{n=0}^\infty$} \label{70625a}

We have already seen in Example \ref{70929a} that frames of the form $\{T^n \varphi\}_{n=0}^\infty$ are quite
sensitive to perturbations. If we restrict ourself to perturb a frame $\{T^n \varphi\}_{n=0}^\infty$ with
elements from a subspace on which $T$ acts as a contraction, a useful stability result can be obtained:

\bpr \label{70601a} Assume that $\{T^n \varphi\}_{n=0}^\infty$ is a frame for some bounded
linear operator $T: \h \to \h$ and some $\varphi \in \h,$
and let $A$ denote a lower frame bound. Assume that $V\subset \h$
is invariant under $T$ and that there exists $\mu\in [0, 1[$ such that
$||T\widetilde{\varphi}||\le \mu\, || \widetilde{\varphi}||$ for all
$\widetilde{\varphi} \in V.$ Then the following hold:

\bei
\item[(i)] $\{T^n(\varphi+\widetilde{\varphi})\}_{n=0}^\infty$ is a frame sequence for all
$\widetilde{\varphi}\in V$.
\item[(ii)]
$\{T^n(\varphi+\widetilde{\varphi})\}_{n=0}^\infty$ is a  frame for all
$\widetilde{\varphi}\in V$   for which
$||\widetilde{\varphi}||<\sqrt{A(1-\mu^2)}.$ \eni \epr

\bp For the proof of (i), by \cite{CH} (or Theorem 22.2.1 in \cite{CB}) it is sufficient to show that the operator
\begin{equation*}
K: \ell^{2}(\Bbb{N}_0)\longrightarrow {\h},\quad K\{c_{n}\}_{n=0}^{\infty}=\sum_{n=0}^{\infty}c_{n}\Big(T^{n}\varphi -T^{n}(\varphi+\widetilde{\varphi})\Big)=\sum_{n=0}^{\infty}c_{n}T^{n}\widetilde{\varphi}
\end{equation*}
is a well-defined and compact operator whenever $\widetilde{\varphi}\in V.$ The assumption on $\widetilde{\varphi}$  implies that $\{T^n \widetilde{\varphi} \}$ is a Bessel sequence, so $K$ is well defined. Now, for
$N\in \mn,$ consider the finite-dimensional operator $K_N: \ell^{2}(\Bbb{N}_0)\longrightarrow {\h}$
given by
$K_{N}\{c_{n}\}_{n=0}^{\infty}:=\sum_{n=0}^{N}c_{n}T^{n}\widetilde{\varphi}.$
Then
\begin{equation*}
\begin{aligned}
\| K-K_{N} \|
&=\sup_{\| \{c_{n}\}\|=1}\| (K-K_{N})\{c_{n}\}_{n=0}^{\infty}\|
\le \| \widetilde{\varphi}\| \left(\sum_{n=N+1}^{\infty}\mu^{2n}\right)^{1/2},
\end{aligned}
\end{equation*} which implies that $\| K-K_{N} \|\to 0$
as $ N \longrightarrow \infty $. Thus the operator $K$ is indeed compact, as desired.

For the proof of (ii),
considering any $\widetilde{\varphi}\in V,$ we have
\bes \sum_{n=0}^\infty \nl  T^n(\varphi+\widetilde{\varphi})-T^n \varphi\nr^2
=  \sum_{n=0}^\infty \nl  T^n\widetilde{\varphi}\nr^2
\le  \sum_{n=0}^\infty \mu^{2n} \nl  \widetilde{\varphi}\nr^2
= \frac{\nl  \widetilde{\varphi}\nr^2}{1-\mu^2}.
\ens Thus, letting $A$ denote a lower frame bound for $\{T^n \varphi\}_{n=0}^\infty$,
we have that
$\sum_{n=0}^\infty \nl  T^n(\varphi+\widetilde{\varphi})-T^n \varphi\nr^2< A$ whenever
$\nl  \widetilde{\varphi}\nr <\sqrt{A(1-\mu^2)}.$ The result now follows from
the perturbation results in \cite{CH,FZ} (or \cite{CB}, page 565). \ep

\bex \label{70601b} Let us return to Example \ref{70614b}.
For any $N\in \mn,$ let $V:= \Span \{e_k\}_{k=1}^N.$ Then
$V$ is invariant under $T,$ and $||T\widetilde{\varphi}|| \le \lambda_N ||\widetilde{\varphi}||$
for all $\widetilde{\varphi}\in V.$ Thus, by Proposition \ref{70601a}
$\{T^n(\varphi+\widetilde{\varphi})\}_{n=0}^\infty$ is a  frame for all
$\widetilde{\varphi}\in V$  with sufficiently small norm.   \ep \enx

\subsection{A no-go result for compact operators} \label{70627c}

It was proved by Aldroubi et al. \cite{A1} that if the operator $T: \h \to \h$ is
compact and self-adjoint, then  $\cup_{j=1}^J\{T^n \varphi_j\}_{n=0}^\infty$ can not
be a frame for $\h$ for any finite collection of vectors
$\varphi_1, \dots, \varphi_J \in \h.$  We will now generalize this result, and show that
the same conclusion holds for all compact operators.

\bpr \label{70212a} Let $\h$ be an infinite-dimensional Hilbert space, assume that
$T: \h \to \h$ is compact, and let $\varphi_1, \dots, \varphi_J \in \h.$ Then $\cup_{j=1}^J\{T^n \varphi_j\}_{n=0}^\infty$ can not
be a frame for $\h.$ \epr

\bp Assume that $\cup_{j=1}^J \{T^n \varphi_j\}_{n=0}^\infty$ is a frame for $\h.$ Then,
by removal of the finitely many elements  $\varphi_1, \dots, \varphi_J,$ the remaining
vectors
$\cup_{j=1}^J \{T^n \varphi_j\}_{n=1}^\infty$ form a frame for the infinite-dimensional
space $V:= \span \{T^n \varphi_j\}_{j=1, \dots, J, n\in \mn}.$ In particular, for
any $f\in V$ there exists some
coefficients $\{c_{n,j}\}_{j=1, \dots, J, n\in \mn}\in \ell^2(\{1, \dots, J\} \times \mn)$
such that
\bes f= \sum_{j=1}^J \sum_{n=1}^\infty c_{n,j} T^n \varphi_j
= T \left(\sum_{j=1}^J \sum_{n=1}^\infty c_{n,j} T^{n-1} \varphi_j\right).\ens
Thus the range of $T$ equals the space $V;$ in particular   $\cra_T$  is closed. In order to arrive at a contradiction,
assume that $T$ is compact.
Since $\h$ is infinite-dimensional
and spanned by the vectors $\varphi_j, T \varphi_j, T^2 \varphi_j, \cdots, \, j=1, \dots, J,$ the range $\cra_T$ is
infinite-dimensional. Consider now the restriction of $T$ to the orthogonal complement
of the kernel of $T,$ i.e.,
$\widetilde{T}: \cn(T)^\bot \to \cra_T.$
Now $\widetilde{T}$ is a bijection; the assumption that $T$ is compact implies that
$\widetilde{T}$ is also compact, but this leads to a contradiction. Let us prove this. In fact, since $\widetilde{T}$ is a bijection between Hilbert spaces,
we know that $\widetilde{T}^{-1}$ is bounded; it follows that for any $f\in \h,$
$||f|| = || \widetilde{T}^{-1}\widetilde{T} f || \le || \widetilde{T}^{-1}||\,
||\widetilde{T} f ||.$
Thus, letting $\etk$ denote an orthonormal basis for $\cn(T)^\bot$ and considering
any $k\neq \ell,$
\bes || \widetilde{T}e_k - \widetilde{T}e_\ell|| \ge \frac1{|| \widetilde{T}^{-1}||}\, ||e_k - e_\ell||
= \frac{\sqrt{2}}{|| \widetilde{T}^{-1}||}.\ens
It follows that $\{\widetilde{T}e_k\}_{k=1}^\infty$ does not have a convergent subsequence,
contradicting the compactness of $\widetilde{T}.$ \ep

\noindent{\bf Ole Christensen, DTU Compute, Technical University of Denmark,
2800 Kgs. Lyngby, Denmark,
Email: ochr@dtu.dk  \\

\noindent Marzieh Hasannasab, DTU Compute, Technical University of Denmark,
2800 Kgs. Lyngby, Denmark,
Email: mhas@dtu.dk \\

\noindent Ehsan Rashidi, Faculty of Mathematical Science, University of Mohaghegh Ardabili, Ardabil, Iran, Email: erashidi@uma.ac.ir }

\end{document}